\documentclass{amsart}

\usepackage{lmodern}
\usepackage[T1]{fontenc}
\usepackage[utf8]{inputenc}

\usepackage{overpic}

\usepackage{amsmath,amssymb,amsfonts,amsthm}
\usepackage[nobysame]{amsrefs}
\usepackage{mathtools}
\usepackage{xspace}

\usepackage{hyperref}
\hypersetup{
  pdfauthor={Zijia Li, Hans-Peter Schröcker, Johannes Siegele, Daren Thimm},
  pdftitle={Quadratic Spinor Polynomials with Infinitely Many Factorizations},
  pdfkeywords={conformal geometric algebra, circular translation, Villarceau
    motion, motion factorization, Dupin cyclide},
  colorlinks=true,
  allcolors=blue,
}

\newcommand{\D}{\mathbb{D}}
\renewcommand{\H}{\mathbb{H}}
\renewcommand{\DH}{\D\H}
\renewcommand{\P}{\mathbb{P}}
\newcommand{\R}{\mathbb{R}}
\newcommand{\Bm}{B_{-}}
\newcommand{\Bp}{B_{+}}
\newcommand{\rv}[1]{\tilde{#1}}

\newcommand{\ci}{\mathrm{i}}
\newcommand{\qi}{\mathbf{i}}
\newcommand{\qj}{\mathbf{j}}
\newcommand{\qk}{\mathbf{k}}
\newcommand{\eps}{\varepsilon}

\newcommand{\CGA}{\ensuremath{\mathrm{CGA}}\xspace}
\newcommand{\CGAp}{\ensuremath{\CGA_{+}}\xspace}
\newcommand{\tp}{\intercal}
\newcommand{\NQ}{\mathcal{N}}
\newcommand{\SV}{\mathcal{S}}

\DeclareMathOperator{\dist}{dist}

\theoremstyle{remark}
\newtheorem{remark}{Remark}

\title{Quadratic Spinor Polynomials with Infinitely Many Factorizations}

\date{\today}

\author{Zijia Li}
\address{KLMM, Academy of Mathematics and Systems Science,
  Chinese Academy of Sciences, Beijing, China}
\email{lizijia@amss.ac.cn}

\author{Hans-Peter Schröcker}
\address{Department of Basic Sciences in Engineering, University of Innsbruck,
  Innsbruck, Austria}
\email{hans-peter.schroecker@uibk.ac.at}

\author{Johannes Siegele}
\address{RICAM Johann Radon Institute for Computational and Applied Mathematics
  Austrian Academy of Sciences, Linz, Austria}
\email{johannes.siegele@oeaw.ac.at}

\author{Daren A. Thimm}
\address{Department of Basic Sciences in Engineering, University of Innsbruck,
  Innsbruck, Austria}
\email{daren.thimm@uibk.ac.at}


\begin{document}

\begin{abstract}
    Spinor polynomials are polynomials with coefficients in the even sub-algebra
  of conformal geometric algebra whose norm polynomial is real. They describe
  rational conformal motions. Factorizations of spinor polynomial corresponds to
  the decomposition of the rational motion into elementary motions. Generic
  spinor polynomials allow for a finite number of factorizations. We present two
  examples of quadratic spinor polynomials that admit infinitely many
  factorizations. One of them, the circular translation, is well-known. The
  other one has only been introduced recently but in a different context. We
  not only compute all factorizations of these conformal motions but also
  interpret them geometrically.
 \end{abstract}

\maketitle

\section{Introduction}
\label{sec:introduction}

\begin{figure}
  \centering
  \includegraphics[]{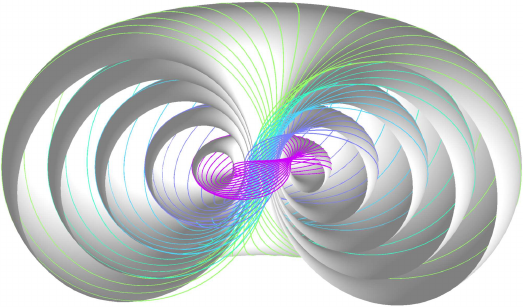}
  \caption{Hopf fibration via Villarceau circles on a hyperbolic family of torus
    surfaces.}
  \label{fig:hopf-fibration}
\end{figure}

In 2019 L.~Dorst presented a conformal motion with rather peculiar geometric
properties \cite{dorst19}:
\begin{itemize}
\item The trajectories of all points are circles.
\item Any point of conformal three-space lies on exactly one circle.
\item The set of all trajectory circles forms the famous Hopf fibration
  \cite{hopf31,zamboj21}
\end{itemize}
The Hopf fibration of space, that is, the trajectories of the motion presented
by Dorst, is visualized in Figure~\ref{fig:hopf-fibration}. In the particular
conformal normal form underlying this article, they can be grouped into families
of Villarceau circles on torus surfaces. Hence, we refer to the motion as
\emph{Villarceau motion.}

The same motion was re-discovered independently by three of the authors of the
present paper during the preparation of \cite{li23} when it attracted our
attention because of curious factorization properties. It can be parametrized by
a quadratic polynomial with coefficients in the conformal geometric algebra
(\CGA) that allows for a two-parametric set of factorizations with linear
factors. This high number combined with its low algebraic degree distinguishes
the conformal Villarceau motion from all other examples of rational motions with
only one exception: The fairly simple planar rigid body motion that translates
the moving object along a circle.

The purpose of this paper is to study and compare algebraic and geometric
properties of factorizations of circular translation and conformal Villarceau
motion. It will turn out that there are some similarities but also important
differences. Ultimately, the conformal Villarceau motion sticks out as a prime
example of a simple yet non-trivial rational motion with an outstanding number
of fascinating geometric, algebraic, and kinematic properties.

This paper continues with a concise introduction to conformal geometric algebra,
factorization theory of spinor polynomials over this algebra, and a kinematic
map for conformal geometry in Section~\ref{sec:preliminaries}. We then proceed
with a study of the factorizations of the circular translation and geometric
interpretations in Section~\ref{sec:circular-translation} before passing on to
the conformal Villarceau motion in Section~\ref{sec:villarceau-motion}. We
recall Dorst's original definition, compute all of its factorizations, and
discuss them in terms of the Euclidean geometry of the individual factors, in
kinematic terms and also as conic section in the projective image space of
conformal kinematics.

\section{Preliminaries}
\label{sec:preliminaries}

This section intends to give a concise introduction to concepts used in the
remainder of this article. For more detailed explanations we refer the reader to
the references mentioned in the text.

Rigid body transformations are generated by an even number of reflections in
planes. Similarly, conformal transformations in space can be defined as
compositions of an even number of reflections (inversions) in spheres or planes
(which usually are considered as “spheres” through the infinity point $\infty$).
The non-linear nature of these transformations suggests that matrices may not
always be the optimal tool to describe them. Moreover, the factorization theory
of rational motions is not naturally formulated in matrix algebras. Hence we use
the framework of conformal geometric algebra (\CGA). Here we will only be
regarding the three-dimensional case.

\subsection{Conformal Geometric Algebra}
\label{sec:cga}

To construct $\CGA$ we first choose an orthonormal basis \(\{e_1, e_2, e_3,
e_+, e_-\} \in \mathbb{R}^{4, 1}\) fulfilling the following conditions
\begin{equation*}
e_1^2=e_2^2=e_3^2=e_+^2=1,\ e_-^2=-1.
\end{equation*}
The multiplication on this algebra is defined to be anticommutative:
\begin{equation*}
  e_i e_j=-e_j e_i \coloneqq e_{ij} \qquad\text{for}\quad i \neq j;\ i,j\in \{1,2,3,+,-\}
\end{equation*}
By linear extension this generates a real algebra, which is called the
\emph{conformal geometric algebra.}

On the \CGA we define an involution called \emph{reversion,} defined by
inverting the order of multiplication of the basis elements. Using the
anticommutativity of the multiplication we can write the reversion in
the following way.
\begin{equation*}
  \tilde e_{i_1, \ldots, i_n}= e_{i_n, \ldots, i_1} = (-1)^{\tfrac{n(n-1)}{2}}e_{i_1, \ldots, i_n}
\end{equation*}
A general element of \CGA and it's reverse can be written as
\begin{equation*}
  a = \sum_{l\subseteq \{ 1, 2, 3, +, -\} } a_l e_l,\quad \tilde a = \sum_{l\subseteq \{ 1, 2, 3, +, -\} } a_l \tilde e_l
\end{equation*}
The cardinality of \(l\) is called the \emph{grade} of the element. This is
defined for the basis vectors \(e_l\) and for elements consisting only of
elements of the same grade. Otherwise we say an element has \emph{mixed grade.}

The basic objects we can represent in \CGA are spheres, planes and points.
Spheres with center \(c\) and radius \(r\) are represented by
\begin{equation}
  \label{eq:1}
  s = c_x e_1 + c_y e_2 + c_z e_3 + \tfrac{1}{2}(1+c_x^2 + c_y^2 +c_z^2-r^2) (e_+ + e_-) - e_+.
\end{equation}
Note that $s\rv{s} = \rv{s}s = r^2$ is real. Points can be represented as
spheres with zero radius. Planes can be seen as a limiting case of a sphere with
radius approaching infinity. Hence the representation of planes looks slightly
different. A plane with unit normal vector \(n\) and oriented distance from the
origin \(d\) is given by
\begin{equation}
  \label{eq:2}
  p = n_x e_1 + n_y e_2 + n_z e_3 + d (e_+ + e_-).
\end{equation}
Throughout this text, we will be using homogeneous coordinates for spheres,
points, and planes. That is, spheres, points, and planes are also represented by
non-zero scalar multiples of \eqref{eq:1} or \eqref{eq:2}.

To be able to study kinematics in \CGA we first have to define how an element
can be acted upon. Any element \(a\) of \CGA (a sphere, point or plane) can be
transformed by any other element \(b\) via the sandwich product given by \(a
\mapsto b a \tilde b\). If \(b\) is a sphere or a plane the action of \(b\)
describes a reflection in $b$. Through this we can construct the group of
conformal displacements, as it is generated by reflection in spheres and planes.
We denote the even sub-algebra of \CGA as \CGAp. It contains all linear
combinations of elements of CGA with even grade and contains the compositions of
an even number of reflections.

Acting with points on other objects is possible. It yields non-invertible maps
which, nonetheless, are of high relevance in kinematics and factorization
theory.

\subsection{Spinor Polynomials and Their Factorization}
\label{sec:factorization-theory}

Consider now a polynomial $C = \sum_{i=0}^n c_it^i$ in the indeterminate $t$ and
with coefficients $c_i \in \CGAp$. Multiplication of polynomials is defined by
the convention that the indeterminate $t$ commutes with all coefficients. Right
evaluation of $C$ at $h \in \CGAp$ is defined as
\begin{equation*}
  C(h) = \sum_{i=0}^n c_ih^i.
\end{equation*}
Because of non-commutativity, this is different from left evaluation
$\sum_{i=0}^n h^ic_i$ (which we will not need in this text).

Denote by $\rv{C}$ the polynomial obtained by conjugating its coefficient $c_0$,
$c_1$, \ldots, $c_n$. For $C$ to describe a conformal motion, it is necessary
that $C\rv{C} = \rv{C}C$ is a non-zero \emph{real polynomial.} If this
condition is met and $C$ is of positive degree, we call $C$ a \emph{spinor
  polynomial} \cite{li23}. The action of spheres/points/planes is extended to
the action of $C$ on spheres/points/planes via the formula
\begin{equation*}
  x \mapsto Cx\rv{C}.
\end{equation*}
If $x$ is a point, then the right-hand side is a polynomial curve in homogeneous
coordinates, that is, a rational curve in Cartesian coordinates. Hence one
speaks of a \emph{rational conformal motion.} Any factorization $C = PQ$ of $C$
with spinor polynomials $P$ and $Q$ corresponds to the decomposition of the
motion into simpler sub-motions.

Most important for applications like in kinematics \cite{hegedus13,gallet16} but
also in discrete differential geometry \cite{hoffmann24} are factorizations with
linear spinor polynomial factors. We only consider factorizations of this type
but, for the sake of brevity, only speak of “factorizations”. A lot is already
known about such factorizations \cite{li18,li23}:
\begin{itemize}
\item For generic $C$ of degree $n$ there exist finitely many factorizations.
  The total number depends on the number of real zeros of the polynomial
  $C\rv{C} \in \R[t]$ and ranges between $n!$ for no real zeros of $C\rv{C}$ and
  $(2n)!/(2^n)$ for the maximum of $2n$ real zeros of $C\rv{C}$.
\item There are examples of spinor polynomials that do not admit any
  factorization.
\item There are examples of spinor polynomials that admit infinitely many
  factorizations.
\item The linear polynomial $t-h$ is a right factor of $C$ if and only if $h$ is
  a right zero of $C$. In other words, $C = C’(t-h)$ for some spinor polynomial
  $C’$, if and only if $C(h) = 0$.
\end{itemize}
The last item shows that linear right factors and right zeros are closely
related.

Examples of spinor polynomials with no or with infinitely many factorizations
are rare. The only two known examples (up to conformal equivalence) of quadratic
spinor polynomials with a two-parametric set of factorizations will be studied
in this article.

\subsection{A Map of Conformal Kinematics}
\label{sec:kinematic-map}

The composition of an even number of reflections in spheres or planes is given
by an algebra element $x \in \CGAp$ of real norm, i.e., $x\rv{x} = \rv{x}x \in
\R$. Since it is only unique up to non-zero scalar multiples, it is naturally
viewed as a point of projective space $\P(\CGAp) = \P^{15}(\R)$. In this way, we
obtain a kinematic map from the group of conformal transformations into
$\P^{15}(\R)$.

Compositions from the left or from the right with fixed conformal
transformations generate a transformation group of $\P^{15}(\R)$ and thus
determines a geometry of that space. Geometrically relevant entities in this
context are:
\begin{itemize}
\item The \emph{Study variety $\SV$,} a projective variety of dimension ten and
  degree twelve that is given by the quadratic equations that encode the
  conditions $x\rv{x} = \rv{x}x \in \R$.
\item The \emph{null quadric $\NQ$,} given by vanishing condition of the real
  part of $x\rv{x}$ (which equals the real part of $\rv{x}x$).
\end{itemize}

Points in the intersection of $\SV$ and $\NQ$ can be thought of as “singular”
conformal displacements. While they do not describe proper conformal
transformations, they are of geometric relevance. They arise quite naturally in
conformal kinematics, for example as conformal scalings (see below) where the
scale factor in the limit goes to zero.

\begin{figure}
  \centering
  \includegraphics[]{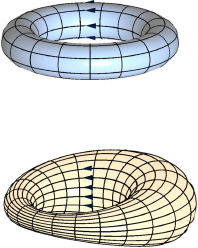}
  \includegraphics[]{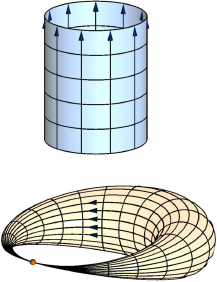}
  \includegraphics[]{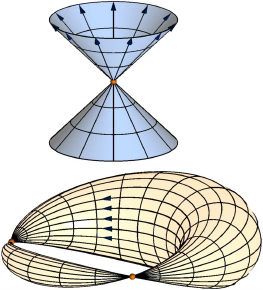}
  \caption{Elementary conformal motions: conformal rotation, translation
    (transversion), and scaling}
  \label{fig:elementary-motions}
\end{figure}

Linear spinor polynomials $at+b$ describe particularly simple rational conformal
motions \cite{dorst16}. Depending on the sign of $h\rv{h}$ with $h = a^{-1}b$,
these elementary motions are classified as follows:
\begin{description}
\item[Conformal rotation ($h\rv{h} > 0$):] A Euclidean rotation around a fixed
  axis and with variable angle or any conformal image thereof.
\item[Conformal translation (transversion; $h\rv{h} = 0$):] A Euclidean
  translation with fixed direction and variable distance or any conformal image
  thereof.
\item[Conformal scaling ($h\rv{h} < 0$):] A uniform scaling from a fixed center
  but with variable scaling factor or any conformal image thereof.
\end{description}
These elementary motions are illustrated in Figure~\ref{fig:elementary-motions}.
We display, in the top row, a typical surface generated by the points of a
circle undergoing a Euclidean rotation, translation, or scaling. In the bottom
row, the trajectory surface (a cyclide of Dupin) of a conformally equivalent
motion is displayed. The arrows indicate point paths.

\section{The Circular Translation}
\label{sec:circular-translation}

In this section we consider a spinor polynomial in the sub-algebra $\langle 1,
\qi, \qj, \qk, \eps \rangle$, where
\begin{equation*}
\qi = -e_{2,3},\quad
  \qj = e_{1,3},\quad
  \qk = -e_{1,2},\quad
  \eps = e_{1,2,3,+} + e_{1,2,3,-}.
\end{equation*}
The generating elements satisfy $\qi^2 = \qj^2 = \qk^2 = \qi\qj\qk = -1$,
$\eps^2 = 0$ and $\qi\eps = \eps\qi$, $\qj\eps = \eps\qj$, $\qk\eps = \eps\qk$,
whence it is the algebra of \emph{dual quaternions} which we denote as $\DH$.
The polynomial we scrutinize is
\begin{equation*}
  C = t^2 + 1 - \eps(\qj t + \qi).
\end{equation*}
Observe that
\begin{equation*}
  C\rv{C} = \rv{C}C = M^2
\end{equation*}
where $M = t^2 + 1$. Hence, $C$ is indeed a spinor polynomial. Since it is
defined over $\DH$, it is even a \emph{motion polynomial} in the sense of
\cite{hegedus13} and it describes a rigid body motion.

In order to write $C$ as product of two linear spinor polynomial factors, we use
a general factorization algorithm \cite{kalkan22,li23}. We use polynomial
division to write $C = QM + R$ where $\deg R < \deg M = 2$. (We will not make
use of the fact that $Q = 1$.) A linear right factor $H_2 = t - h_2$ of $C$ is
necessarily a right factor of $M$ and hence also of the linear remainder
polynomial $R = r_1t + r_0$. Moreover, it should satisfy the spinor polynomial
conditions $H_2\rv{H}_2$, $\rv{H}_2{H_2} \in \mathbb{R}[t]$ which boil down to
$h_2\rv{h}_2$, $\rv{h}_2h_2$, $h_2 + \rv{h}_2 \in \R$, c.f. \cite{kalkan22}. For
general spinor polynomials, these conditions are satisfied precisely by $h_2 =
-r_1^{-1}r_0$. However, in our case, we have $r_1 = -\eps \qj$ which is not
invertible. Existence of a right factor thus hinges on the existence of $h_2 \in
\DH$ subject to the conditions:
\begin{equation}
  \label{eq:3}
  \begin{aligned}
    r_1h_2 + r_0 = 0 &\quad \text{($t-h_2$ is right zero of $R$)} \\
    h_2^2 + 1 = 0 &\quad \text{($h_2$ is a zero of $M$)} \\
    h_2\rv{h}_2, \rv{h}_2h_2, h_2 + \rv{h}_2 \in \R &\quad \text{($t-h_2$ is a spinor polynomial)}
  \end{aligned}
\end{equation}
Using the 16 yet unknown coefficients of $h_2$ as variables, we can convert
these into a system of 39 algebraic equations. Among them, 13 are of degree one
and 26 are of degree two. This system of equations allows for a straightforward
solution. It turns out that all factorizations are over the dual quaternions
$\DH$:
\begin{equation}
  \label{eq:4}
  C = (t - \qk - \eps ((1-\mu)\qj - \lambda\qi))
      (t + \qk - \eps(\lambda\qi + \mu\qj)),
  \quad
  (\lambda, \mu) \in \R^2.
\end{equation}

\begin{remark}
  The factorizations \eqref{eq:4} of the motion polynomial $C$ are well-known
  \cite{gallet16}. We have gained the additional insight that there are not more
  factorizations over $\CGAp$ than there are over~$\DH$.
\end{remark}

Several geometric interpretations of \eqref{eq:4} are conceivable. To begin
with, we may view the factors as parametric surfaces
\begin{equation*}
  h_1(\lambda,\mu) = \qk - \lambda\eps\qi + (1-\mu)\eps\qj,
  \quad
  h_2(\lambda,\mu) = -\qk + \lambda\eps\qi + \mu\eps\qj
\end{equation*}
in the affine spaces parallel to the vectors $\eps\qi$ and $\eps\qj$ and through
$\qk$ and $-\qk$, respectively. \emph{We see that the map $h_2(\lambda,\mu)
  \mapsto h_1(\lambda,\mu)$ is the composition of the reflection in the point
  $-\qk$ followed by the translation by the vector $\eps\qj + 2\qk$.}

\begin{figure}
  \centering
  \includegraphics[page=3]{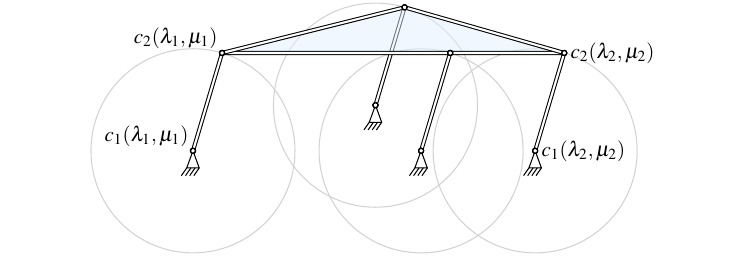}
  \caption{Circular translation and parallelogram linkage}
  \label{fig:parallelogram-linkage}
\end{figure}

A second interpretation of \eqref{eq:4} is in terms of kinematics of the
underlying motion. Each factor describes a rotation around an axis parallel to
$\qk$. Thus, $C$ parametrizes a planar rigid body motion. The respective
rotation centers $c_1(\lambda,\mu)$, $c_2(\lambda,\mu)$ of $t-h_1$ and $t-h_2$ are the points
\begin{equation*}
  c_1(\lambda,\mu) = (1-\mu, \lambda)^\tp,\quad
  c_2(\lambda,\mu) = (-\mu, \lambda)^\tp.
\end{equation*}
The point $c_2(\lambda,\mu)$ actually rotates around $c_1(\lambda,\mu)$. Its trajectory is a circle
parametrized as
\begin{equation*}
  c_2(\lambda,\mu,t) = \frac{1}{1+t^2}(-\mu(1+t^2)+2, \lambda(1+t^2) + 2t)^\tp.
\end{equation*}
Denote by $\dist(a,b)$ the distance between two points $a$ and $b$. It is
straightforward to deduce the following statements:
\begin{itemize}
\item For any $(\lambda,\mu) \in \R^2$ we have $\dist(c_1(\lambda,\mu),
  c_2(\lambda,\mu,t)) = 1$.
\item For any two $(\lambda_1,\mu_1)$, $(\lambda_2,\mu_2) \in \R^2$ we have
  $\dist^2(c_1(\lambda_1,\mu_1), c_1(\lambda_2,\mu_2)) =
  \dist^2(c_2(\lambda_1,\mu_1,t), c_2(\lambda_2,\mu_2,t)) = (\lambda_1-\lambda_2)^2 + (\mu_1-\mu_2)^2$.
\end{itemize}
\emph{We infer that $C$ describes the coupler motion of a parallelogram linkage
  where input and output crank are both of length $1$}
(Figure~\ref{fig:parallelogram-linkage}). Each factorization of $C$ gives rise
to one additional leg of length~$1$. The coupler motion itself is the
translation along a circle.

For generic cases, a geometric relation between factorizations of spinor
polynomials and the geometry in the kinematic image space space $\P^{15}(\R)$ is
well-understood. The interpretation of non-generic spinor polynomials is less
clear. In the example of a circular translation, we can restrict ourselves to
the projective space $\P^3(\R)$ over the vector space and sub-algebra generated
by $1$, $\eps\qi$, $\eps\qj$, and $\qk$. We confine ourselves to observing that
the rational curve parametrized by $C$ is, indeed, rather special.

The transformation group generated by left and right multiplications (coordinate
changes in fixed and moving frame) turns $\P^3$ into a quasi elliptic space
\cite[Section~3.2.2]{klawitter15}. Using homogeneous coordinates $x = x_0 +
x_1\eps\qi + x_2\eps\qj + x_3\qk$, its invariant figure consists of the two
planes with equations
\begin{equation*}
  \omega\colon x_0 + \ci x_3 = 0,
  \quad
  \overline{\omega}\colon x_0 - \ci x_3 = 0
\end{equation*}
(the intersection of the image space of planar kinematics with the null quadric
$\NQ$) and the two points $[n_1]$, $[n_2]$ where
\begin{equation*}
  n_1 = \eps\qi + \ci\eps\qj,
  \quad
  n_2 = \eps\qi - \ci\eps\qj
\end{equation*}
(c.f. \cite[Chapter~11,~\S2]{bottema90}).\footnote{Here, we freely extend
  projective space and algebra to complex coefficients. Note that the complex
  unit $\ci$ needs to be distinguished from the quaternion unit $\qi$.}

The polynomial $C$ of \eqref{eq:4} parametrizes a rational quadratic curve (a
conic section) in this quasi elliptic space. It lies in the plane $x_3 = 0$ of
pure translations. Moreover, $C(\ci) = n_1$ and $C(-\ci) = n_2$ so that the
conic section contains the absolute points of quasi elliptic space. We may hence
address it as a \emph{quasi elliptic circle.}

\section{The Conformal Villarceau Motion}
\label{sec:villarceau-motion}

The conformal Villarceau motion was introduced by L.~Dorst in \cite{dorst19}. In
contrast to the circular translation of the previous section, it is not a
rigid-body motion. Some of its curious geometric properties have been briefly
covered in Section~\ref{sec:introduction}

We recall the parametric equation of the Villarceau motion of Dorst in
\cite{dorst19}. With
\begin{equation*}
  \Bm \coloneqq e_{12}
  \quad\text{and}\quad
  \Bp \coloneqq e_{3+}
\end{equation*}
it is given by
\begin{equation*}
  C = \exp(-\Bm\tfrac{\varphi}{2})\exp(-\Bp\tfrac{\varphi}{2}).
\end{equation*}
Expanding this in a Taylor series, simplifying using, Using $\Bm^2 = \Bp^2 =
-1$, and separating even and odd parts, this can be written as
\begin{equation*}
  C = (\cos\varphi + \Bm \sin\varphi)(\cos\varphi + \Bp \sin\varphi).
\end{equation*}
Substituting $t = -\cot{\varphi}$ we arrive at the polynomial
\begin{equation*}
  C = (t - \Bm)(t - \Bp) = t^2 - t(e_{12} + e_{3+}) + e_{123+}.
\end{equation*}
It admits one factorization with linear spinor polynomial factors by
construction. We wish to find all of these factorizations. In doing so, we
follow the general procedure that has already been outlined in
Section~\ref{sec:circular-translation}. Once more, the norm polynomial is
\begin{equation*}
  C\rv{C} = \rv{C}C = M^2
\end{equation*}
where $M = (t^2+1)^2$. Using polynomial division we write $C = QM + R$ where $R
= r_1t + r_0$ and
\begin{equation*}
  r_1 = -e_{12}-e_{3+},\quad
  r_0 = e_{123+}-1.
\end{equation*}
Again, $r_1$ is not invertible so that we have to solve the system of
equations~\eqref{eq:3}. Using the 16 yet unknown coefficients of $h_2$ as
variables, we can convert these into a system of 43 algebraic equations. Among
them, 17 are of degree one and 26 are of degree two. Solving the 17 linear
equations results in
\begin{equation}
  \label{eq:5}
  h_2 = e_{12} + s_x x + s_y y + s_z z
\end{equation}
where
\begin{gather}
  \label{eq:6}
  s_x = 2(e_{1+}-e_{23}),\quad
  s_y = 2(e_{2+}+e_{13}),\quad
  s_z = 2(e_{3+}-e_{12}),\\
  \label{eq:7}
  \text{and}\quad
  x^2 + y^2 + (z - \tfrac{1}{4})^2 - \tfrac{1}{16} = 0.
\end{gather}
Since the vectors $e_{1+}-e_{23}$, $e_{2+}+e_{13}$, $e_{3+}-e_{12}$ are pairwise
perpendicular and of equal length we may say that $h_2$ lies on a sphere given
by \eqref{eq:5}--\eqref{eq:7}. If $S(u,v)$ is some parametrization of this
sphere, we have
\begin{equation}
  \label{eq:8}
  h_2(u,v) = m + \tfrac{1}{4}S(u,v)
\end{equation}
where $m = \tfrac{1}{2}(e_{12} + e_{3+})$. Using polynomial division, we find $C
= (t - h_1(u,v))(t - h_2(u,v))$ where
\begin{equation}
  \label{eq:9}
  h_1(u,v) = m - \tfrac{1}{4}S(u,v).
\end{equation}
Thus, we can say: \emph{The factors $t-h_1(u,v)$ and $t-h_2(u,v)$ in the
  factorizations of the Villarceau motion are parameterized by the points of a
  sphere. The map $h_1(u,v) \mapsto h_2(u,v)$ is the reflection in the sphere
  center $m$. This implies commutativity of $h_1(u,v)$ and $h_2(u,v)$ which can,
  of course, also be confirmed by a straightforward calculation.}

Elementary conformal motions are more difficult to grasp intuitively than mere
rotations or translations. Hence, a kinematic explanation of the infinitely many
factorizations of $C = (t - h_1(u,v))(t - h_2(u,v))$ is a challenging task. Lets
first look at the individual elementary motions $H_1 = s - h_1(u,v)$ and $H_2 =
t - h_2(u,v)$, parametrized by different parameters $s$ and $t$. In this sense,
their product creates a two-parametric rational motion. It is readily confirmed
that
\begin{equation}
  \label{eq:10}
  H_1H_2 = H_2H_1.
\end{equation}
For any point $x$, the trajectory surface $D_x = H_1H_2x\rv{H}_2\rv{H}_1$ is a
cyclide of Dupin. This can be seen by
\begin{enumerate}
\item verifying that the parameter lines are circles and showing that
\item the second fundamental form of $D_x$ is diagonal.
\end{enumerate}
This implies that the curvature lines of $D_x$ are circles which is a
characteristic property of cyclides of Dupin.

\begin{figure}
  \centering
\includegraphics[page=2,trim=50 0 30 0,clip,scale=0.48]{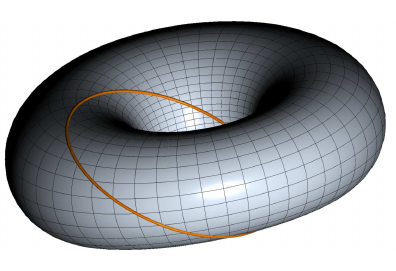}\includegraphics[page=3,trim=30 0 50 0,clip,scale=0.48]{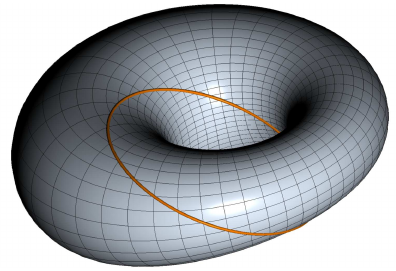}\includegraphics[page=4,trim=25 0 60 0,clip,scale=0.48]{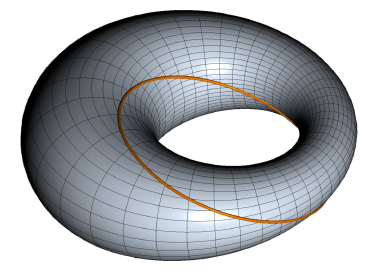}\includegraphics[page=5,trim=53 0 70 0,clip,scale=0.48]{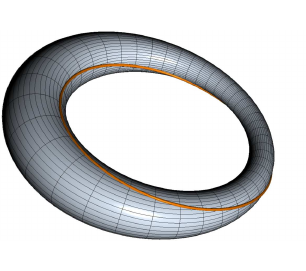}
  \caption{Different Dupin cyclides with the same Villarceau circle}
  \label{fig:dupin-cyclides}
\end{figure}

The parameter line circles are point trajectories under the elementary motions
given by $H_1$ and $H_2$, respectively. The trajectory of point $x$ under the
one-parametric rational motion $C$ is obtained as diagonal surface curve for $s
= t$. In a conformal transformation of $D_x$ to a torus, it is mapped to a
Villarceau circle. Hence, we follow \cite{dorst19} and address it as
\emph{Villarceau circle} on a Dupin cyclide. Different values of $u$ and $v$
provide different Dupin cyclides but \emph{the same Villarceau circle} on them.
This is illustrated in Figure~\ref{fig:dupin-cyclides}.

The interpretation of the conformal Villarceau motion $C$ in kinematic image
space is more involved than in case of the circular translation and we cannot
hope to resolve all mysteries. The geometry of kinematic image space is not
among the classical non-Euclidean geometries and a complete system of geometric
invariants is not known. Therefore, we confine ourselves to highlight some
peculiarities.

The polynomial $C$ parametrizes a rational curve of degree two in $\P^{15}(\R)$.
It intersects the null quadric $\NQ$ in only two points $[n_1]$, $[n_2]$ where
\begin{equation*}
  n_1 = C(\ci) = e_{123+} - 1 - \ci(e_{12} + e_{3+}),
  \quad
  n_2 = C(-\ci) = e_{123+} - 1 + \ci(e_{12} + e_{3+}).
\end{equation*}
It is known \cite{scharler21} that factorizability is related the connecting
lines of intersection points of $C$ and $\NQ$. In our case, these are the conic
tangents in $[n_1]$ and $[n_2]$ as well as their connecting line. More
precisely, the linear remainder polynomial $R$ in the expression $C = QM + R$
parametrizes one such line. We already observed that its leading coefficient
$r_1$ is not invertible. A direct computation shows that \emph{none of the
  points on the connecting line of $[n_1]$ and $[n_2]$ or on the tangents of $C$
  in $[n_1]$ and $[n_2]$ are invertible.} This explains why conventional
factorization attempts (with only finitely many factorizations) fail in case of
$C$ and it emphasizes the importance of non-invertible elements of \CGAp. Not
much seems to be known about their geometry as points of~$\P(\CGAp)$.

\section{Conclusion}
\label{sec:conclusion}

We have investigated the two known examples of quadratic spinor polynomials with
a two parameteric set of factorizations. The underlying motions are the fairly
simple translation along a circle and the Villarceau motion of \cite{dorst19}.
Both motions have curious algebraic and geometric properties and exhibit some
similarities but also differences. Most notably, the factors of the Villarceau
motion are parametrized by the points of a sphere and not a plane. They
correspond in the reflection in the sphere center and do commute, c.f.
\eqref{eq:10}. While the circular translation can be viewed as a quasi elliptic
circle in kinematic image space, the analogous interpretation of the Villarceau
motion is less obvious. For the circular translation the point on the secant
connecting its null points are non-invertible, for the Villarceau motion this is
the case along the secant but also along the tangents. The circular translation
can be created mechanically by a parallelogram linkage with factorizations
corresponding to possible cranks. For the factorizations of the conformal
Villarceau motion and its trajectories we gave an interpretation in terms of
Dupin cyclides that share a Villarceau circle.

A more systematic study of the conformal Villarceau motion is planned for the
future. The question for low degree spinor polynomials with many factorizations
apart from circular translation and conformal Villarceau motion remains wide
open.

\section*{Acknowledgment}

We gratefully acknowledge support by the Austrian Science Fund (FWF) P~33397-N
(Rotor Polynomials: Algebra and Geometry of Conformal Motions).
 
\bibliographystyle{spmpsci}

\end{document}